\documentclass{amsart}
\usepackage[english]{babel}
\usepackage[latin1]{inputenc}
\usepackage[dvips,final]{graphics}
\usepackage{amsmath,amsfonts,amssymb,amsthm,amscd,array,stmaryrd,mathrsfs, mathdots, epigraph}
\usepackage{pstricks}
 \usepackage[all]{xy}
 \usepackage{url}
\usepackage{multirow, blkarray}
\usepackage{booktabs}
\usepackage{textcomp}
 \usepackage[final]{epsfig}
\newtheorem{lemma}{Lemma}[subsection]
\newtheorem{proposition}[lemma]{Proposition}
\newtheorem{remark}[lemma]{Remark}
\newtheorem{example}[lemma]{Example}
\newtheorem{theorem}{Theorem}
\newtheorem{definition}[lemma]{Definition}
\newtheorem{corollary}[lemma]{Corollary}

\newcommand{\eps}{{\varepsilon}}

\newcommand{\om}{{\omega}}

\newcommand{\proofend}{$\Box$\bigskip} 

\newcommand{\Diff}{{\mathrm {Diff}}}
\newcommand{\PSL}{{\mathrm {PSL}}}
\newcommand{\SL}{{\mathrm {SL}}}
\newcommand{\Id}{{\mathrm {Id}}}

\newcommand{\R}{{\mathbb R}}
\newcommand{\Z}{{\mathbb Z}}
\newcommand{\RP}{{\mathbb {RP}}}  

\def\proof{\paragraph{Proof.}}

\begin{document}

\title{Coxeter's frieze patterns and discretization of the Virasoro orbit }

\author[V. Ovsienko]{Valentin Ovsienko}

\author[S. Tabachnikov]{Serge Tabachnikov}

\address{
Valentin Ovsienko,
CNRS,
Laboratoire de Math\'ematiques, 
Universit\'e de Reims-Champagne-Ardenne, 
FR~3399 CNRS, F-51687, Reims, France}

\address{
Serge Tabachnikov,
Pennsylvania State University,
Department of Mathematics,
University Park, PA 16802, USA,
and 
ICERM, Brown University, Box1995,
Providence, RI 02912, USA
}

\email{
ovsienko@math.univ-lyon1.fr,
tabachni@math.psu.edu
}

\date{}


\maketitle

\begin{abstract}
We show that the space of classical Coxeter's frieze patterns
can be viewed as a discrete version of a coadjoint orbit of the
Virasoro algebra.
The canonical (cluster) (pre)symplectic form on the space
of frieze patterns is a discretization of the Kirillov symplectic form. We relate a continuous version of frieze patterns to conformal metrics of constant curvature in dimension 2.
\end{abstract}

\maketitle
\date{}

\section{Introduction and main results} \label{Intro}

The goal of this note is to relate two different subjects:

\begin{enumerate}
\item
the {\it Virasoro algebra} and the related infinite-dimensional symplectic
manifold $\Diff_+(S^1)/\PSL_2(\R)$;

\item
the space of {\it Coxeter's frieze patterns} viewed as a cluster manifold.

\end{enumerate}

The problem of discretization of the Virasoro algebra
is very well known and has been studied by different authors,
and different discretizations were suggested~\cite{FT,FRS};
see also~\cite{MS}.
The main motivation for this study is application to
integrable systems, such as the Korteweg - de Vries (KdV) equation.
Several discrete versions of the KdV were proposed.
Most of the discrete versions of the Virasoro algebra
consist of discretization of the corresponding linear Poisson structure
on its dual space.

We will describe a discretization procedure that relates the subject
to combinatorics and cluster algebra.
We will be interested in the infinite-dimensional homogeneous space $\Diff_+(S^1)/\PSL_2(\R)$
equipped with (a $1$-parameter family of) Kirillov's symplectic structures.  
This symplectic space is often regarded as a coadjoint orbit of the Virasoro algebra~\cite{K1,K2},
or, in other words, a symplectic leaf of the linear Poisson structure.
This is a more geometric way to understand the Virasoro-related Poisson structure.
We will obtain a finite-dimensional discretization of $\Diff_+(S^1)/\PSL_2(\R)$.

We do not consider integrable systems in this paper,
but believe that our discretization procedure can be applied to KdV
and should be related to such discrete integrable systems as the pentagram map;
see~\cite{OST} and references therein.

The discrete objects that we consider are the classical
{\it Coxeter (or Conway-Coxeter) frieze patterns}~\cite{Cox,CoCo}.
This notion was invented in the early 1970's but became widely known quite recently
due to its close relation to the {\it cluster algebra};
see~\cite{CaCh,ARS}.
In particular, it was shown in \cite{MGOT,MGOST} 
that frieze patterns are closely related with the moduli space of polygons 
in the projective line and with second order linear difference equations. 
As every cluster manifold, the space of Coxeter's friezes has a canonical
(pre)symplectic structure.
We will prove the following.

\begin{theorem}
\label{MainOne}
The space $\Diff_+(S^1)/\PSL_2(\R)$ equipped with Kirillov's symplectic form 
is a continuous limit of the space of Coxeter's frieze patterns equipped with the
cluster (pre)symplectic form.
\end{theorem}

Following the ideas of Conway and Coxeter, we identify frieze patterns
with linear recurrence equations:
\begin{equation}
\label{Req}
V_{i+1} = c_i V_i - V_{i-1},
\end{equation}
where the ``potential'' $(c_i)$ is an $n$-periodic sequence of (real or complex) numbers, 
and where the sequence $(V_i)$ is unknown, i.e., a ``solution''.
Furthermore, we will impose
the following condition: the potential is $n$-periodic, and all the solutions are $n$-antiperiodic:
\begin{equation}
\label{Per}
c_{i+n}=c_i,
\qquad
V_{i+n}=-V_i.
\end{equation}
The space of equations~(\ref{Req}) satisfying the condition~(\ref{Per})
is an algebraic variety of dimension $n-3$.
If $n$ is odd, then this algebraic variety is isomorphic to the classical
moduli space $\mathcal{M}_{0,n}$; see~\cite{MGOST}
for  details.

Note that equation~(\ref{Req}) is nothing other than the classical 
{\it discrete Hill equation} (also called Sturm-Liouville or Schr\"odinger, equation).
Relation of this equation to a discrete version of the Virasoro algebra
is very natural and appeared in all the above cited works on the subject.
However, the notion of frieze pattern and cluster algebra were not considered.
We will show that this approach provides additional combinatorial tools.
Let us also emphasize the fact that the (anti)periodicity condition~(\ref{Per})
seems to be the only natural way to obtain a finite-dimensional
space of equations~(\ref{Req}) approximating the space $\Diff_+(S^1)/\PSL_2(\R)$.

We also show how to describe the continuous limit of
Coxeter's frieze patterns in terms of solutions of the
classical Liouville equation.
Solutions of this equation are interpreted in terms of 
projective differential geometry.
We believe that geometric and combinatorial viewpoints
complement each other and lead to a better understanding of both
parts of the story.

We end the introduction by an open question that concerns
a natural generalization of Coxeter's friezes called $2$-friezes; see~\cite{MGOT}.
The space of $2$-friezes is an algebraic variety of dimension $2n-8$.
It is related to linear recurrence equations
of order $3$:
$$
V_{i+3}=a_iV_{i+2} - b_i V_{i+1} +V_{i},
$$
with $n$-periodic solutions.
This space also carries a structure of cluster manifold,
and therefore has a canonical (pre)symplectic structure.
This is the space on which the pentagram map acts.
It would be natural to expect that the canonical cluster symplectic structure
is related to the Gelfand-Dickey bracket.
However, the pentagram map does not preserve the canonical symplectic structure.
It would be very interesting to understand the situation in that case.

\section{The space of Coxeter's friezes} \label{CCox}

In this section, we recall the classical notion of a Coxeter frieze pattern.
We introduce local coordinate systems and identify this space
with the moduli space $\mathcal{M}_{0,n}$.

\subsection{Closed frieze patterns}
We start with the definition of  classical Coxeter's frieze patterns.

\begin{definition}
{\rm
\begin{enumerate}
\item[(a)]
A frieze pattern~\cite{Cox} is an infinite array of numbers 
$$
\begin{array}{cccccccccccc}
&\cdots&&0&&0&&0&&0&&\cdots\\[4pt]
\cdots&&1&&1&&1&&1&&\cdots\\[4pt]
&\cdots&&c_{i}&&c_{i+1}&&c_{i+2}&&c_{i+3}&&\cdots\\[4pt]
&& \cdots&& \cdots&& \cdots&& \cdots&&
\end{array}
$$
where the entries propagate downward, and  the entries of each next row are  
determined by the previous two rows via the frieze rule.
For each elementary ``diamond''
\begin{equation}
\label{diamond1}
\begin{array}{ccc}
&b&\\[2pt]
a&&d\\[2pt]
&c&
\end{array}
\end{equation}
one has
\begin{equation}
\label{RulEq}
ad-bc=1.
\end{equation}
For instance, the entries in the next row of the above frieze are
$c_ic_{i+1}-1$. 

\item[(b)]
A frieze pattern is called {\it closed} if a row of $1$'s appears again:
$$
\begin{array}{ccccccccccc}
\cdots&&1&&1&&1&&1&&\cdots\\[4pt]
&c_i&&c_{i+1}&&c_{i+2}&&c_{i+3}&&c_{i+4}\\[4pt]
&& \cdots&& \cdots&& \cdots&& \cdots&&\\[4pt]
\cdots&&1&&1&&1&&1&&\cdots\\[4pt]
&0&&0&&0&&0&&\cdots\\[4pt]
\cdots&&-1&&-1&&-1&&-1&&\cdots
\end{array}
$$
By definition, this lower row of $1$'s is  followed by a row of~$0$'s, and then by a row of $-1$'s. 
One can extend the array vertically so that each diagonal, 
in either of the two directions, is {\it anti-periodic}.

\item[(c)]
The {\it width} $w$ of a closed frieze pattern is the number of non-trivial rows between the rows of~$1$'s. 
\end{enumerate}
}
\end{definition}

\begin{example}
\label{CEx}
{\rm
A generic Coxeter frieze pattern of width~$2$ is as follows:
$$
 \begin{array}{ccccccccccc}
\cdots&&1&& 1&&1&&\cdots
 \\[4pt]
&a_1&&\frac{a_2+1}{a_1}&&\frac{a_1+1}{a_2}&&a_2&&
 \\[4pt]
\cdots&&a_2&&\frac{a_1+a_2+1}{a_1a_2}&&a_1&&\cdots
 \\[4pt]
&1&&1&&1&&1&&
\end{array}
$$
for some $a_1,a_2\not=0$.
(Note that we omitted the first and the last rows of $0$'s.)}
\end{example}

The following facts are well-known~\cite{Cox,CoCo}; see also~\cite{MGOST}.

\begin{enumerate}
\item
A closed frieze pattern is horizontally periodic with period
$$
n=w+3,
$$ that is, $c_{i+n}=c_i$. 
\item
A frieze pattern with the first line $(c_i)$ is closed if and only if
the equation~(\ref{Req}) satisfies  condition~(\ref{Per}).
\item
A closed frieze pattern has ``glide symmetry'' whose second iteration 
is the horizontal parallel translation distance~$n$. 
\end{enumerate}

The name ``frieze pattern" is due to the glide symmetry.

\subsection{Local coordinates}

Every Coxeter's frieze pattern of width $w$ is uniquely defined by its
(South-East) diagonal:
\begin{equation}
\label{CanDiag}
 \begin{array}{ccccccccccc}
1&&1&& 1&&1
 \\[4pt]
&a_1&&\cdots&&&&
 \\[4pt]
&&a_2&&\cdots&&
 \\[4pt]
 &&&\ddots&&&\\[4pt]
 &&\cdots&&a_w&&\cdots\\[4pt]
&1&&1&&1&&1&&
\end{array}
\end{equation}
for some $a_1,a_2,\ldots,a_w\not=0$ (see Example~\ref{CEx}).
Therefore, $(a_1,a_2,\ldots,a_w)$ is a local coordinate system
on the space of friezes.

Consider a different choice of the diagonal,
or, more generally, consider an arbitrary {\it zigzag}
\begin{equation}
\label{ZZCo} 
\begin{array}{ccccccccccc}
1&&1&& 1&&1
 \\[4pt]
&a'_1&&\cdots&&&&
 \\[4pt]
&&a'_2&&\cdots&&
 \\[4pt]
&a'_3&&\cdots&&
 \\[4pt]
a'_4&&\cdots&&
 \\[4pt]
 &\ddots&&&\\[4pt]
 &&a'_w&&\cdots\\[4pt]
&1&&1&&1&&1&&
\end{array}
\end{equation}
Again the frieze pattern is uniquely defined by  
$a'_1,a'_2,\ldots,a'_w\not=0$,
so that we obtain a different coordinate system $(a'_1,a'_2,\ldots,a'_w)$.

It turns out that the coordinate changes between these coordinate
systems can be understood as {\it mutations} in the cluster algebra
of type $A_w$; see~\cite{MGOT} (the Appendix). 
The space of all Coxeter's friezes is therefore a cluster manifold.

\subsection{Relation to the moduli space $\mathcal{M}_{0,n}$}
The following statement is proved in~\cite{MGOT,MGOST}.

\begin{proposition}
\label{PropM} 
If $n$ is odd, then the space of Coxeter's friezes is isomorphic 
to the moduli space $\mathcal{M}_{0,n}$.
\end{proposition}

\proof
We briefly describe the main construction; see~\cite{MGOT,MGOST}
for the details.

Consider Coxeter's frieze pattern given by a diagonal~(\ref{CanDiag}).
Take the neighboring (upper) diagonal and write the two diagonals together
to obtain $n$ vectors in~$\R^2$:
\begin{equation}
\label{CProj} 
\left(\begin{array}{c}
0\\[4pt]
1
\end{array}\right),
\quad
\left(\begin{array}{c}
1\\[4pt]
a_1
\end{array}\right),
\quad
\left(\begin{array}{c}
\frac{a_2+1}{a_1}\\[4pt]
a_2
\end{array}\right),
\quad\ldots,\quad
\left(\begin{array}{c}
1\\[4pt]
0
\end{array}\right).
\end{equation}
Note that these vectors $V_0,\ldots,V_{n-1}$ form a fundamental solution
of the equation~(\ref{Req}) corresponding to the frieze~(\ref{CanDiag}).

Projecting this $n$-gon to $\RP^1$, one obtains a point in $\mathcal{M}_{0,n}$,
and this projection is a one-to-one correspondence.
\proofend

\begin{remark}
\label{ImpRem} 
{\rm
The above construction allows us to explain the geometric meaning of the entries of the frieze.
The elements $(a_1,\ldots,a_w)$ of the diagonal
can be obtained as the vector products:
$$
a_i=\left[V_{n-1},V_i\right],
$$
for $i=1,\ldots,w$.
}
\end{remark}

\section{Continuous Coxeter's friezes} \label{fri}

In this section, we introduce our main notion
of a continuous frieze pattern that can be obtained as a continuous limit 
of classical Coxeter friezes.
We show that the continuous limit of a Coxeter's frieze can be understood
as a solution of the Liouville equation with special boundary conditions.
The space of these solutions is identified with $\Diff_+(S^1)/\PSL_2(\R)$.

\subsection{Projective curves, Hill's equations and the space $\Diff_+(S^1)/\PSL_2(\R)$}

Let $\Diff_+(S^1)$ be the group of orientation preserving diffeomorphisms
of the circle $S^1\simeq\RP^1$.
The homogeneous space $\Diff_+(S^1)/\PSL_2(\R)$ is
one of the most interesting infinite-dimensional manifolds
 in geometry and mathematical physics.
Its study was initiated by Kirillov~\cite{K1,K2}.
We give here two different realizations of this space.
The results of this section are well-known.

\begin{definition}
{\rm
We call a {\it (simple) projective curve} an orientation preserving diffeomorphism 
$$
\gamma:\R/T\Z \to \RP^1,
$$ 
that is, a parameterization of the projective line by $[0,T)$. 
The {\it projective equivalence class} of $\gamma$ consists of the diffeomorphisms $\varphi\circ\gamma$ where
$\varphi:\RP^1\to\RP^1$ is a projective transformation, i.e., $\varphi\in\PSL_2(\R)$. 
}
\end{definition}

The space of projective equivalence classes of curves is isomorphic to
$\Diff_+(S^1)/\PSL_2(\R)$, the $\Diff_+(S^1)$-action on the curves being given by
\begin{equation} 
\label{ActC}
f:\gamma\mapsto\gamma\circ{}f^{-1}
\qquad
\hbox{where}
\qquad
f\in\Diff_+(S^1).
\end{equation}
This space can also be identified with projective structures on $\RP^1$ with monodromy $-\Id$;
see~\cite{K1,K2} and also~\cite{OT}.

\begin{definition}
{\rm
\begin{enumerate}
\item[(a)]
A {\it Hill equaltion} is a $2$nd order linear differential equation of the form
\begin{equation} 
\label{Hill}
2c\,y''(x)+k(x)y(x)=0
\end{equation}
with $T$-periodic potential $k(x)$, here $c\in\R$ is an arbitrary constant.
We will always assume that Hill's equation has monodromy $-\Id$,
i.e., all the solutions of (\ref{Hill}) are $T$-anti-periodic:
$$
y(x+T)=-y(x).
$$

\item[(b)]
A Hill equation (\ref{Hill}) with $T$-anti-periodic solutions 
is called {\it non-oscillating} if every
solution has exactly $1$ zero on $[0,T)$.
\end{enumerate}
}
\end{definition}
An example of a non-oscillating Hill equation with $T$-anti-periodic solutions is the equation with the constant potential
$k(x)\equiv{2c\pi^2}/{T^2}$.

\begin{remark}
{\rm
Recall that
the diagonals of a frieze pattern are solutions to a linear difference equation of $2$nd order~(\ref{Req})
where $c_i$ are the terms in the first non-trivial row of the frieze pattern. 
Hill's equation is a continuous analog of this difference equation, 
and its potential, $k(x)$ is a continuous analog of the sequence $(c_i)$.
}
\end{remark}

The following statement can be found in~\cite{K1}; see also~\cite{OT}
for a detailed discussion. 
We give a proof for the sake of completeness.

\begin{proposition}
\label{KirProp}
The space $\Diff_+(S^1)/\PSL_2(\R)$
can be identified with the space of non-oscillating
of Hill's equations.
\end{proposition}

\proof
 The space of solutions of Hill's equation (\ref{Hill}) is two-dimensional, and the Wronski determinant of any pair of solutions is a constant that we normalize to be equal to 1. 
 Thus a choice of two solutions determines a curve 
 $\Gamma(x) \subset \R^2$ such that $[\Gamma,\Gamma'] =1$ (the bracket denotes the determinant formed by two vectors). 
 This curve $\Gamma(x)$ is well defined up to the action of~$\SL_2(\R)$,
it is antiperiodic, i.e., $\Gamma(x+T)=-\Gamma(x)$.
Furthermore, the curve is simple if and only if the equation is non-oscillating. 
 
 Conversely, given a projective curve $\gamma$,
 we can lift $\gamma$ to a parameterized curve $\Gamma(x) \subset \R^2$. 
The closure condition on $\gamma$ implies that $\Gamma(x+T)=-\Gamma(x)$. 
The lift is not unique: one can always multiply $\Gamma(x)$ by a non-vanishing function~$\lambda(x)$. 
We fix the lift by the condition
\begin{equation} 
\label{deter}
[\Gamma(x),\Gamma'(x)]=1\quad {\rm for\ all}\ x.
\end{equation}
Projectively equivalent curves $\gamma$ correspond to $\SL(2,\R)$-equivalent curves $\Gamma$, see \cite{OT}.
Differentiating (\ref{deter}), implies that the vectors $\Gamma$ and $\Gamma''$ are proportional,
i.e., $\Gamma$ satisfies Hill's equaltion
$$
\Gamma''(x)=k(x)\Gamma(x)
$$
with $T$-periodic potential $k(x)$ and monodromy $-\Id$.
\proofend

One can recover the potential $k(x)$ of Hill's equation from the curve $\Gamma$: if $f(x)$ is the ratio of two coordinates of the curve $\Gamma(x)$ then
$$
k=c\, S(f)
\qquad 
\hbox{where}
\qquad S(f)=\frac{f'''}{f'}-\frac{3}{2}\left(\frac{f''}{f'}\right)^2
$$
is the classical {\it Schwarzian derivative}.
This formula is quite old and should probably be attributed to Lagrange; see~\cite{What}.
A more contemporary way to express the same observation is
the following formula of $\Diff_+(S^1)$-action on the space of Hill's equations:
\begin{equation} 
\label{AdStar}
f:k(x)\mapsto{}k(f^{-1}(x))\left({f^{-1}}'\right)^2+ c\, S(f^{-1})(x).
\end{equation}
Note that this more complicated formula precisely corresponds to the action~(\ref{ActC}).

\subsection{Definition of continuous friezes}\label{ConLimS}
Let us now describe the procedure of {\it continuous limit}
of a frieze pattern.

A natural labeling of the entries of a frieze pattern is according to the scheme:
\begin{equation} \label{diamond}
\begin{array}{ccc}
&v_{i,j}&\\[4pt]
v_{i,j-1}&&v_{i+1,j}\\[4pt]
&v_{i+1,j-1}&
\end{array}
\end{equation}

A continuous analog of a frieze pattern is a twice differentiable function of two variables $F(x,y)$ 
satisfying an analog of the frieze rule. Namely, replace (\ref{diamond}) by
$$
\begin{array}{ccc}
&F(x,y+\eps)&\\[4pt]
F(x,y)&&F(x+\eps,y+\eps)\\[4pt]
&F(x+\eps,y)&
\end{array}
$$
and expand in $\eps$ up to the 2nd order. Then the frieze rule yields
$$
\eps^2(F F_{xy} - F_{x} F_{y})=1,
$$
where $F_{x}$ and $F_{y}$ denote the partial derivatives with respect to $x$ and $y$ respectively. 
Thus, we are led to the following. 

\begin{definition}
{\rm
\begin{enumerate}
\item[(a)]
A {\it continuous frieze pattern} is a function $F(x,y)$ satisfying the partial differential equation
\begin{equation} 
\label{Liou}
F(x,y) F_{xy}(x,y)-F_{x}(x,y) F_{y}(x,y)=1.
\end{equation}

\item[(b)]
A closed continuous frieze patterns is
a function $F(x,y)$ satisfying~(\ref{Liou}) and the following conditions:
\begin{equation} 
\label{bdry}
\begin{array}{rl}
F(x,x) =0,\ F_{y}(x,x)=1 &{\rm for\ all}\ x;\\[6pt]
F(x,y)>0 & {\rm for}\ x<y<x+T;\\ [6pt]
F(x+T,y)=F(x,y+T)=-F(x,y)& {\rm for\ all}\ x,y.
\end{array}
\end{equation}
\end{enumerate}
}
\end{definition}

\begin{remark}
{\rm
Equation~(\ref{Liou}) is the classical {\it Liouville equation} on the function $\ln F$, see, e.g.,~\cite{DFN}.
The first two conditions in (\ref{bdry}) are analogs of having a row of $0$'s, 
followed by a  row of $1$'s. 
The last condition is an analog of anti-periodicity of the diagonals, $T$ being the period. 
Given the first condition, the  second one is equivalent to $F_{x}(x,x) = -1$.}
\end{remark}

\subsection{Continuous friezes from projective curves}\label{ConSSec}
Let us describe a simple geometric construction that provides 
all the solutions to~(\ref{Liou}) satisfying~(\ref{bdry}).
We shall show that a projective equivalence class
of a projective curve determines a continuous frieze pattern,
and that every continuous frieze pattern can be obtained from
a projective curve.

Consider a projective curve $\gamma$ and its canonical lift $\Gamma$ to $\R^2$
satisfying  condition~(\ref{deter}).

\begin{theorem} 
\label{curve}
(i)
The function 
\begin{equation} 
\label{Constr}
F(x,y)=[\Gamma(x),\Gamma(y)]
\end{equation}
is a closed continuous frieze pattern.

(ii)
Conversely, every closed continuous frieze pattern is of the form~(\ref{Constr}) for some curve~$\gamma$.
\end{theorem}

\proof 
Part (i). One has:
$$
F_{x}=[\Gamma'(x),\ \Gamma(y)], \qquad
F_{y}=[\Gamma(x),\ \Gamma'(y)], \qquad
F_{xy}=[\Gamma'(x),\Gamma'(y)].
$$
The Ptolemy (or Pl\"ucker) relation for the determinants made by the vectors $\Gamma(x), \Gamma(y), \Gamma'(x), \Gamma'(y)$ implies
$$
[\Gamma (x),\Gamma(y)]  [\Gamma'(x),\Gamma'(y)] - 
 [\Gamma'(x),\Gamma(y)]  [\Gamma(x),\Gamma'(y)] = 
  [\Gamma(x),\Gamma'(x)]  [\Gamma(y),\Gamma'(y)],
$$
and (\ref{Liou}) follows.

The first and the last of the boundary conditions (\ref{bdry}) obviously hold, and the second condition coincides with (\ref{deter}). The positivity follows from the fact that $\Gamma(x)$ induces an embedding of the interval $(0,T)$ to $\RP^1$.

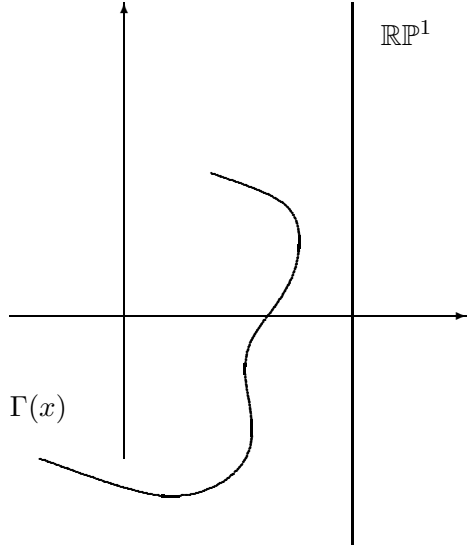
\begin{center}
\setlength{\unitlength}{3144sp}%
\begingroup\makeatletter\ifx\SetFigFont\undefined%
\gdef\SetFigFont#1#2#3#4#5{%
  \reset@font\fontsize{#1}{#2pt}%
  \fontfamily{#3}\fontseries{#4}\fontshape{#5}%
  \selectfont}%
\fi\endgroup%
\begin{picture}(3624,4319)(4714,-6158)
\put(4726,-5123){\makebox(0,0)[lb]{\smash{{\SetFigFont{11}{16.8}{\rmdefault}{\mddefault}{\updefault}{\color[rgb]{0,0,0}$\Gamma(x)$}%
}}}}
\thinlines
{\color[rgb]{0,0,0}\put(4726,-4336){\vector( 1, 0){3600}}
}%
\thicklines
\multiput(4958,-5464)(9.00000,-3.00000){2}{\makebox(6.3500,9.5250){\SetFigFont{7}{8.4}{\rmdefault}{\mddefault}{\updefault}.}}
\multiput(4967,-5467)(14.10000,-4.70000){2}{\makebox(6.3500,9.5250){\SetFigFont{7}{8.4}{\rmdefault}{\mddefault}{\updefault}.}}
\multiput(4981,-5472)(9.39655,-3.75862){3}{\makebox(6.3500,9.5250){\SetFigFont{7}{8.4}{\rmdefault}{\mddefault}{\updefault}.}}
\multiput(5000,-5479)(8.33333,-3.33333){4}{\makebox(6.3500,9.5250){\SetFigFont{7}{8.4}{\rmdefault}{\mddefault}{\updefault}.}}
\multiput(5025,-5489)(7.80000,-2.60000){5}{\makebox(6.3500,9.5250){\SetFigFont{7}{8.4}{\rmdefault}{\mddefault}{\updefault}.}}
\multiput(5056,-5500)(9.07500,-3.02500){5}{\makebox(6.3500,9.5250){\SetFigFont{7}{8.4}{\rmdefault}{\mddefault}{\updefault}.}}
\multiput(5092,-5513)(8.28000,-2.76000){6}{\makebox(6.3500,9.5250){\SetFigFont{7}{8.4}{\rmdefault}{\mddefault}{\updefault}.}}
\multiput(5133,-5528)(7.55000,-2.51667){7}{\makebox(6.3500,9.5250){\SetFigFont{7}{8.4}{\rmdefault}{\mddefault}{\updefault}.}}
\multiput(5178,-5544)(8.05000,-2.68333){7}{\makebox(6.3500,9.5250){\SetFigFont{7}{8.4}{\rmdefault}{\mddefault}{\updefault}.}}
\multiput(5226,-5561)(8.40000,-2.80000){7}{\makebox(6.3500,9.5250){\SetFigFont{7}{8.4}{\rmdefault}{\mddefault}{\updefault}.}}
\multiput(5276,-5579)(8.55000,-2.85000){7}{\makebox(6.3500,9.5250){\SetFigFont{7}{8.4}{\rmdefault}{\mddefault}{\updefault}.}}
\multiput(5327,-5597)(8.50000,-2.83333){7}{\makebox(6.3500,9.5250){\SetFigFont{7}{8.4}{\rmdefault}{\mddefault}{\updefault}.}}
\multiput(5378,-5614)(8.35000,-2.78333){7}{\makebox(6.3500,9.5250){\SetFigFont{7}{8.4}{\rmdefault}{\mddefault}{\updefault}.}}
\multiput(5428,-5631)(8.35000,-2.78333){7}{\makebox(6.3500,9.5250){\SetFigFont{7}{8.4}{\rmdefault}{\mddefault}{\updefault}.}}
\multiput(5478,-5648)(7.80000,-2.60000){7}{\makebox(6.3500,9.5250){\SetFigFont{7}{8.4}{\rmdefault}{\mddefault}{\updefault}.}}
\multiput(5525,-5663)(7.60000,-2.53333){7}{\makebox(6.3500,9.5250){\SetFigFont{7}{8.4}{\rmdefault}{\mddefault}{\updefault}.}}
\multiput(5571,-5677)(8.70000,-2.90000){6}{\makebox(6.3500,9.5250){\SetFigFont{7}{8.4}{\rmdefault}{\mddefault}{\updefault}.}}
\multiput(5615,-5690)(8.47058,-2.11764){6}{\makebox(6.3500,9.5250){\SetFigFont{7}{8.4}{\rmdefault}{\mddefault}{\updefault}.}}
\multiput(5657,-5702)(7.85882,-1.96470){6}{\makebox(6.3500,9.5250){\SetFigFont{7}{8.4}{\rmdefault}{\mddefault}{\updefault}.}}
\multiput(5696,-5713)(9.29412,-2.32353){5}{\makebox(6.3500,9.5250){\SetFigFont{7}{8.4}{\rmdefault}{\mddefault}{\updefault}.}}
\multiput(5733,-5723)(8.70588,-2.17647){5}{\makebox(6.3500,9.5250){\SetFigFont{7}{8.4}{\rmdefault}{\mddefault}{\updefault}.}}
\multiput(5768,-5731)(8.47060,-2.11765){5}{\makebox(6.3500,9.5250){\SetFigFont{7}{8.4}{\rmdefault}{\mddefault}{\updefault}.}}
\multiput(5802,-5739)(7.98078,-1.59615){5}{\makebox(6.3500,9.5250){\SetFigFont{7}{8.4}{\rmdefault}{\mddefault}{\updefault}.}}
\multiput(5834,-5745)(10.00000,-1.66667){4}{\makebox(6.3500,9.5250){\SetFigFont{7}{8.4}{\rmdefault}{\mddefault}{\updefault}.}}
\multiput(5864,-5750)(9.62163,-1.60361){4}{\makebox(6.3500,9.5250){\SetFigFont{7}{8.4}{\rmdefault}{\mddefault}{\updefault}.}}
\multiput(5893,-5754)(9.24323,-1.54054){4}{\makebox(6.3500,9.5250){\SetFigFont{7}{8.4}{\rmdefault}{\mddefault}{\updefault}.}}
\put(5921,-5757){\line( 1, 0){ 27}}
\put(5948,-5759){\line( 1, 0){ 27}}
\put(5975,-5761){\line( 1, 0){ 26}}
\put(6001,-5761){\line( 1, 0){ 29}}
\put(6030,-5760){\line( 1, 0){ 29}}
\multiput(6059,-5759)(9.56757,1.59459){4}{\makebox(6.3500,9.5250){\SetFigFont{7}{8.4}{\rmdefault}{\mddefault}{\updefault}.}}
\multiput(6088,-5756)(9.24323,1.54054){4}{\makebox(6.3500,9.5250){\SetFigFont{7}{8.4}{\rmdefault}{\mddefault}{\updefault}.}}
\multiput(6116,-5753)(9.35137,1.55856){4}{\makebox(6.3500,9.5250){\SetFigFont{7}{8.4}{\rmdefault}{\mddefault}{\updefault}.}}
\multiput(6144,-5748)(9.35897,1.87179){4}{\makebox(6.3500,9.5250){\SetFigFont{7}{8.4}{\rmdefault}{\mddefault}{\updefault}.}}
\multiput(6172,-5742)(9.35897,1.87179){4}{\makebox(6.3500,9.5250){\SetFigFont{7}{8.4}{\rmdefault}{\mddefault}{\updefault}.}}
\multiput(6200,-5736)(8.90000,2.96667){4}{\makebox(6.3500,9.5250){\SetFigFont{7}{8.4}{\rmdefault}{\mddefault}{\updefault}.}}
\multiput(6227,-5728)(9.30000,3.10000){4}{\makebox(6.3500,9.5250){\SetFigFont{7}{8.4}{\rmdefault}{\mddefault}{\updefault}.}}
\multiput(6255,-5719)(8.62070,3.44828){4}{\makebox(6.3500,9.5250){\SetFigFont{7}{8.4}{\rmdefault}{\mddefault}{\updefault}.}}
\multiput(6281,-5709)(9.02300,3.60920){4}{\makebox(6.3500,9.5250){\SetFigFont{7}{8.4}{\rmdefault}{\mddefault}{\updefault}.}}
\multiput(6308,-5698)(8.26667,4.13333){4}{\makebox(6.3500,9.5250){\SetFigFont{7}{8.4}{\rmdefault}{\mddefault}{\updefault}.}}
\multiput(6333,-5686)(8.26667,4.13333){4}{\makebox(6.3500,9.5250){\SetFigFont{7}{8.4}{\rmdefault}{\mddefault}{\updefault}.}}
\multiput(6358,-5674)(7.94117,4.76470){4}{\makebox(6.3500,9.5250){\SetFigFont{7}{8.4}{\rmdefault}{\mddefault}{\updefault}.}}
\multiput(6382,-5660)(7.94117,4.76470){4}{\makebox(6.3500,9.5250){\SetFigFont{7}{8.4}{\rmdefault}{\mddefault}{\updefault}.}}
\multiput(6406,-5646)(7.38460,4.92307){4}{\makebox(6.3500,9.5250){\SetFigFont{7}{8.4}{\rmdefault}{\mddefault}{\updefault}.}}
\multiput(6428,-5631)(7.04000,5.28000){4}{\makebox(6.3500,9.5250){\SetFigFont{7}{8.4}{\rmdefault}{\mddefault}{\updefault}.}}
\multiput(6449,-5615)(6.66667,5.33333){4}{\makebox(6.3500,9.5250){\SetFigFont{7}{8.4}{\rmdefault}{\mddefault}{\updefault}.}}
\multiput(6469,-5599)(6.52460,5.43717){4}{\makebox(6.3500,9.5250){\SetFigFont{7}{8.4}{\rmdefault}{\mddefault}{\updefault}.}}
\multiput(6488,-5582)(8.25000,8.25000){3}{\makebox(6.3500,9.5250){\SetFigFont{7}{8.4}{\rmdefault}{\mddefault}{\updefault}.}}
\multiput(6505,-5566)(5.83333,5.83333){4}{\makebox(6.3500,9.5250){\SetFigFont{7}{8.4}{\rmdefault}{\mddefault}{\updefault}.}}
\multiput(6522,-5548)(7.25410,8.70492){3}{\makebox(6.3500,9.5250){\SetFigFont{7}{8.4}{\rmdefault}{\mddefault}{\updefault}.}}
\multiput(6537,-5531)(6.66000,8.88000){3}{\makebox(6.3500,9.5250){\SetFigFont{7}{8.4}{\rmdefault}{\mddefault}{\updefault}.}}
\multiput(6550,-5513)(6.42000,8.56000){3}{\makebox(6.3500,9.5250){\SetFigFont{7}{8.4}{\rmdefault}{\mddefault}{\updefault}.}}
\multiput(6563,-5496)(5.42645,9.04408){3}{\makebox(6.3500,9.5250){\SetFigFont{7}{8.4}{\rmdefault}{\mddefault}{\updefault}.}}
\multiput(6574,-5478)(5.29410,8.82350){3}{\makebox(6.3500,9.5250){\SetFigFont{7}{8.4}{\rmdefault}{\mddefault}{\updefault}.}}
\multiput(6584,-5460)(4.50000,9.00000){3}{\makebox(6.3500,9.5250){\SetFigFont{7}{8.4}{\rmdefault}{\mddefault}{\updefault}.}}
\multiput(6593,-5442)(3.65515,9.13787){3}{\makebox(6.3500,9.5250){\SetFigFont{7}{8.4}{\rmdefault}{\mddefault}{\updefault}.}}
\multiput(6601,-5424)(3.70000,11.10000){3}{\makebox(6.3500,9.5250){\SetFigFont{7}{8.4}{\rmdefault}{\mddefault}{\updefault}.}}
\multiput(6609,-5402)(3.65000,10.95000){3}{\makebox(6.3500,9.5250){\SetFigFont{7}{8.4}{\rmdefault}{\mddefault}{\updefault}.}}
\multiput(6616,-5380)(2.30770,11.53850){3}{\makebox(6.3500,9.5250){\SetFigFont{7}{8.4}{\rmdefault}{\mddefault}{\updefault}.}}
\multiput(6621,-5357)(1.33333,8.00000){4}{\makebox(6.3500,9.5250){\SetFigFont{7}{8.4}{\rmdefault}{\mddefault}{\updefault}.}}
\multiput(6625,-5333)(1.32433,7.94600){4}{\makebox(6.3500,9.5250){\SetFigFont{7}{8.4}{\rmdefault}{\mddefault}{\updefault}.}}
\multiput(6628,-5309)(1.32433,7.94600){4}{\makebox(6.3500,9.5250){\SetFigFont{7}{8.4}{\rmdefault}{\mddefault}{\updefault}.}}
\put(6631,-5285){\line( 0, 1){ 26}}
\put(6632,-5259){\line( 0, 1){ 26}}
\put(6632,-5233){\line( 0, 1){ 26}}
\put(6632,-5207){\line( 0, 1){ 27}}
\put(6630,-5180){\line( 0, 1){ 27}}
\multiput(6628,-5153)(-1.48650,8.91900){4}{\makebox(6.3500,9.5250){\SetFigFont{7}{8.4}{\rmdefault}{\mddefault}{\updefault}.}}
\multiput(6625,-5126)(-1.48650,8.91900){4}{\makebox(6.3500,9.5250){\SetFigFont{7}{8.4}{\rmdefault}{\mddefault}{\updefault}.}}
\multiput(6622,-5099)(-1.48650,8.91900){4}{\makebox(6.3500,9.5250){\SetFigFont{7}{8.4}{\rmdefault}{\mddefault}{\updefault}.}}
\multiput(6619,-5072)(-1.44143,8.64860){4}{\makebox(6.3500,9.5250){\SetFigFont{7}{8.4}{\rmdefault}{\mddefault}{\updefault}.}}
\multiput(6615,-5046)(-1.44143,8.64860){4}{\makebox(6.3500,9.5250){\SetFigFont{7}{8.4}{\rmdefault}{\mddefault}{\updefault}.}}
\multiput(6611,-5020)(-1.38740,8.32440){4}{\makebox(6.3500,9.5250){\SetFigFont{7}{8.4}{\rmdefault}{\mddefault}{\updefault}.}}
\multiput(6607,-4995)(-1.60257,8.01283){4}{\makebox(6.3500,9.5250){\SetFigFont{7}{8.4}{\rmdefault}{\mddefault}{\updefault}.}}
\multiput(6602,-4971)(-1.91890,11.51340){3}{\makebox(6.3500,9.5250){\SetFigFont{7}{8.4}{\rmdefault}{\mddefault}{\updefault}.}}
\multiput(6598,-4948)(-1.91890,11.51340){3}{\makebox(6.3500,9.5250){\SetFigFont{7}{8.4}{\rmdefault}{\mddefault}{\updefault}.}}
\multiput(6594,-4925)(-1.82430,10.94580){3}{\makebox(6.3500,9.5250){\SetFigFont{7}{8.4}{\rmdefault}{\mddefault}{\updefault}.}}
\multiput(6591,-4903)(-2.09615,10.48075){3}{\makebox(6.3500,9.5250){\SetFigFont{7}{8.4}{\rmdefault}{\mddefault}{\updefault}.}}
\multiput(6587,-4882)(-1.64865,9.89190){3}{\makebox(6.3500,9.5250){\SetFigFont{7}{8.4}{\rmdefault}{\mddefault}{\updefault}.}}
\multiput(6585,-4862)(-1.66215,9.97290){3}{\makebox(6.3500,9.5250){\SetFigFont{7}{8.4}{\rmdefault}{\mddefault}{\updefault}.}}
\multiput(6582,-4842)(-1.89190,11.35140){3}{\makebox(6.3500,9.5250){\SetFigFont{7}{8.4}{\rmdefault}{\mddefault}{\updefault}.}}
\multiput(6580,-4819)(-1.81080,10.86480){3}{\makebox(6.3500,9.5250){\SetFigFont{7}{8.4}{\rmdefault}{\mddefault}{\updefault}.}}
\put(6578,-4797){\line( 0, 1){ 22}}
\put(6578,-4775){\line( 0, 1){ 22}}
\put(6578,-4753){\line( 0, 1){ 22}}
\put(6579,-4731){\line( 0, 1){ 22}}
\multiput(6580,-4709)(1.74325,10.45950){3}{\makebox(6.3500,9.5250){\SetFigFont{7}{8.4}{\rmdefault}{\mddefault}{\updefault}.}}
\multiput(6583,-4688)(1.74325,10.45950){3}{\makebox(6.3500,9.5250){\SetFigFont{7}{8.4}{\rmdefault}{\mddefault}{\updefault}.}}
\multiput(6586,-4667)(2.09615,10.48075){3}{\makebox(6.3500,9.5250){\SetFigFont{7}{8.4}{\rmdefault}{\mddefault}{\updefault}.}}
\multiput(6590,-4646)(2.61765,10.47060){3}{\makebox(6.3500,9.5250){\SetFigFont{7}{8.4}{\rmdefault}{\mddefault}{\updefault}.}}
\multiput(6595,-4625)(3.30000,9.90000){3}{\makebox(6.3500,9.5250){\SetFigFont{7}{8.4}{\rmdefault}{\mddefault}{\updefault}.}}
\multiput(6601,-4605)(3.30000,9.90000){3}{\makebox(6.3500,9.5250){\SetFigFont{7}{8.4}{\rmdefault}{\mddefault}{\updefault}.}}
\multiput(6607,-4585)(3.75860,9.39650){3}{\makebox(6.3500,9.5250){\SetFigFont{7}{8.4}{\rmdefault}{\mddefault}{\updefault}.}}
\multiput(6614,-4566)(3.75860,9.39650){3}{\makebox(6.3500,9.5250){\SetFigFont{7}{8.4}{\rmdefault}{\mddefault}{\updefault}.}}
\multiput(6621,-4547)(3.65515,9.13787){3}{\makebox(6.3500,9.5250){\SetFigFont{7}{8.4}{\rmdefault}{\mddefault}{\updefault}.}}
\multiput(6629,-4529)(4.30000,8.60000){3}{\makebox(6.3500,9.5250){\SetFigFont{7}{8.4}{\rmdefault}{\mddefault}{\updefault}.}}
\multiput(6638,-4512)(4.30000,8.60000){3}{\makebox(6.3500,9.5250){\SetFigFont{7}{8.4}{\rmdefault}{\mddefault}{\updefault}.}}
\multiput(6647,-4495)(4.72060,7.86767){3}{\makebox(6.3500,9.5250){\SetFigFont{7}{8.4}{\rmdefault}{\mddefault}{\updefault}.}}
\multiput(6656,-4479)(5.00000,7.50000){3}{\makebox(6.3500,9.5250){\SetFigFont{7}{8.4}{\rmdefault}{\mddefault}{\updefault}.}}
\multiput(6666,-4464)(4.85295,8.08825){3}{\makebox(6.3500,9.5250){\SetFigFont{7}{8.4}{\rmdefault}{\mddefault}{\updefault}.}}
\multiput(6676,-4448)(5.00000,7.50000){3}{\makebox(6.3500,9.5250){\SetFigFont{7}{8.4}{\rmdefault}{\mddefault}{\updefault}.}}
\multiput(6686,-4433)(5.38460,8.07690){3}{\makebox(6.3500,9.5250){\SetFigFont{7}{8.4}{\rmdefault}{\mddefault}{\updefault}.}}
\multiput(6697,-4417)(6.00000,8.00000){3}{\makebox(6.3500,9.5250){\SetFigFont{7}{8.4}{\rmdefault}{\mddefault}{\updefault}.}}
\multiput(6709,-4401)(6.42000,8.56000){3}{\makebox(6.3500,9.5250){\SetFigFont{7}{8.4}{\rmdefault}{\mddefault}{\updefault}.}}
\multiput(6722,-4384)(6.42000,8.56000){3}{\makebox(6.3500,9.5250){\SetFigFont{7}{8.4}{\rmdefault}{\mddefault}{\updefault}.}}
\multiput(6735,-4367)(6.66000,8.88000){3}{\makebox(6.3500,9.5250){\SetFigFont{7}{8.4}{\rmdefault}{\mddefault}{\updefault}.}}
\multiput(6748,-4349)(7.50000,9.00000){3}{\makebox(6.3500,9.5250){\SetFigFont{7}{8.4}{\rmdefault}{\mddefault}{\updefault}.}}
\multiput(6763,-4331)(7.08000,9.44000){3}{\makebox(6.3500,9.5250){\SetFigFont{7}{8.4}{\rmdefault}{\mddefault}{\updefault}.}}
\multiput(6777,-4312)(5.00000,6.66667){4}{\makebox(6.3500,9.5250){\SetFigFont{7}{8.4}{\rmdefault}{\mddefault}{\updefault}.}}
\multiput(6792,-4292)(5.00000,6.66667){4}{\makebox(6.3500,9.5250){\SetFigFont{7}{8.4}{\rmdefault}{\mddefault}{\updefault}.}}
\multiput(6807,-4272)(5.16000,6.88000){4}{\makebox(6.3500,9.5250){\SetFigFont{7}{8.4}{\rmdefault}{\mddefault}{\updefault}.}}
\multiput(6822,-4251)(5.16000,6.88000){4}{\makebox(6.3500,9.5250){\SetFigFont{7}{8.4}{\rmdefault}{\mddefault}{\updefault}.}}
\multiput(6837,-4230)(4.66667,7.00000){4}{\makebox(6.3500,9.5250){\SetFigFont{7}{8.4}{\rmdefault}{\mddefault}{\updefault}.}}
\multiput(6851,-4209)(4.92307,7.38460){4}{\makebox(6.3500,9.5250){\SetFigFont{7}{8.4}{\rmdefault}{\mddefault}{\updefault}.}}
\multiput(6866,-4187)(4.23530,7.05883){4}{\makebox(6.3500,9.5250){\SetFigFont{7}{8.4}{\rmdefault}{\mddefault}{\updefault}.}}
\multiput(6879,-4166)(4.38237,7.30394){4}{\makebox(6.3500,9.5250){\SetFigFont{7}{8.4}{\rmdefault}{\mddefault}{\updefault}.}}
\multiput(6892,-4144)(4.38237,7.30394){4}{\makebox(6.3500,9.5250){\SetFigFont{7}{8.4}{\rmdefault}{\mddefault}{\updefault}.}}
\multiput(6905,-4122)(3.73333,7.46667){4}{\makebox(6.3500,9.5250){\SetFigFont{7}{8.4}{\rmdefault}{\mddefault}{\updefault}.}}
\multiput(6917,-4100)(3.66667,7.33333){4}{\makebox(6.3500,9.5250){\SetFigFont{7}{8.4}{\rmdefault}{\mddefault}{\updefault}.}}
\multiput(6928,-4078)(3.10343,7.75858){4}{\makebox(6.3500,9.5250){\SetFigFont{7}{8.4}{\rmdefault}{\mddefault}{\updefault}.}}
\multiput(6938,-4055)(4.90000,9.80000){3}{\makebox(6.3500,9.5250){\SetFigFont{7}{8.4}{\rmdefault}{\mddefault}{\updefault}.}}
\multiput(6947,-4035)(4.00000,10.00000){3}{\makebox(6.3500,9.5250){\SetFigFont{7}{8.4}{\rmdefault}{\mddefault}{\updefault}.}}
\multiput(6955,-4015)(4.17240,10.43100){3}{\makebox(6.3500,9.5250){\SetFigFont{7}{8.4}{\rmdefault}{\mddefault}{\updefault}.}}
\multiput(6963,-3994)(3.65000,10.95000){3}{\makebox(6.3500,9.5250){\SetFigFont{7}{8.4}{\rmdefault}{\mddefault}{\updefault}.}}
\multiput(6970,-3972)(3.65000,10.95000){3}{\makebox(6.3500,9.5250){\SetFigFont{7}{8.4}{\rmdefault}{\mddefault}{\updefault}.}}
\multiput(6977,-3950)(2.00000,8.00000){4}{\makebox(6.3500,9.5250){\SetFigFont{7}{8.4}{\rmdefault}{\mddefault}{\updefault}.}}
\multiput(6983,-3926)(2.88235,11.52940){3}{\makebox(6.3500,9.5250){\SetFigFont{7}{8.4}{\rmdefault}{\mddefault}{\updefault}.}}
\multiput(6989,-3903)(1.66667,8.33333){4}{\makebox(6.3500,9.5250){\SetFigFont{7}{8.4}{\rmdefault}{\mddefault}{\updefault}.}}
\multiput(6994,-3878)(1.38740,8.32440){4}{\makebox(6.3500,9.5250){\SetFigFont{7}{8.4}{\rmdefault}{\mddefault}{\updefault}.}}
\multiput(6998,-3853)(1.44143,8.64860){4}{\makebox(6.3500,9.5250){\SetFigFont{7}{8.4}{\rmdefault}{\mddefault}{\updefault}.}}
\put(7002,-3827){\line( 0, 1){ 25}}
\put(7004,-3802){\line( 0, 1){ 26}}
\put(7006,-3776){\line( 0, 1){ 26}}
\put(7007,-3750){\line( 0, 1){ 25}}
\put(7006,-3725){\line( 0, 1){ 25}}
\multiput(7005,-3700)(-1.32433,7.94600){4}{\makebox(6.3500,9.5250){\SetFigFont{7}{8.4}{\rmdefault}{\mddefault}{\updefault}.}}
\multiput(7002,-3676)(-1.33333,8.00000){4}{\makebox(6.3500,9.5250){\SetFigFont{7}{8.4}{\rmdefault}{\mddefault}{\updefault}.}}
\multiput(6998,-3652)(-2.30770,11.53850){3}{\makebox(6.3500,9.5250){\SetFigFont{7}{8.4}{\rmdefault}{\mddefault}{\updefault}.}}
\multiput(6993,-3629)(-2.76470,11.05880){3}{\makebox(6.3500,9.5250){\SetFigFont{7}{8.4}{\rmdefault}{\mddefault}{\updefault}.}}
\multiput(6987,-3607)(-3.50000,10.50000){3}{\makebox(6.3500,9.5250){\SetFigFont{7}{8.4}{\rmdefault}{\mddefault}{\updefault}.}}
\multiput(6980,-3586)(-4.90000,9.80000){3}{\makebox(6.3500,9.5250){\SetFigFont{7}{8.4}{\rmdefault}{\mddefault}{\updefault}.}}
\multiput(6971,-3566)(-4.80000,9.60000){3}{\makebox(6.3500,9.5250){\SetFigFont{7}{8.4}{\rmdefault}{\mddefault}{\updefault}.}}
\multiput(6961,-3547)(-5.42645,9.04408){3}{\makebox(6.3500,9.5250){\SetFigFont{7}{8.4}{\rmdefault}{\mddefault}{\updefault}.}}
\multiput(6950,-3529)(-6.00000,9.00000){3}{\makebox(6.3500,9.5250){\SetFigFont{7}{8.4}{\rmdefault}{\mddefault}{\updefault}.}}
\multiput(6938,-3511)(-6.00000,8.00000){3}{\makebox(6.3500,9.5250){\SetFigFont{7}{8.4}{\rmdefault}{\mddefault}{\updefault}.}}
\multiput(6926,-3495)(-7.25000,7.25000){3}{\makebox(6.3500,9.5250){\SetFigFont{7}{8.4}{\rmdefault}{\mddefault}{\updefault}.}}
\multiput(6912,-3480)(-8.16395,6.80329){3}{\makebox(6.3500,9.5250){\SetFigFont{7}{8.4}{\rmdefault}{\mddefault}{\updefault}.}}
\multiput(6896,-3466)(-8.45900,7.04917){3}{\makebox(6.3500,9.5250){\SetFigFont{7}{8.4}{\rmdefault}{\mddefault}{\updefault}.}}
\multiput(6879,-3452)(-9.44000,7.08000){3}{\makebox(6.3500,9.5250){\SetFigFont{7}{8.4}{\rmdefault}{\mddefault}{\updefault}.}}
\multiput(6860,-3438)(-7.23077,4.82051){4}{\makebox(6.3500,9.5250){\SetFigFont{7}{8.4}{\rmdefault}{\mddefault}{\updefault}.}}
\multiput(6838,-3424)(-7.69607,4.61764){4}{\makebox(6.3500,9.5250){\SetFigFont{7}{8.4}{\rmdefault}{\mddefault}{\updefault}.}}
\multiput(6815,-3410)(-8.80000,4.40000){4}{\makebox(6.3500,9.5250){\SetFigFont{7}{8.4}{\rmdefault}{\mddefault}{\updefault}.}}
\multiput(6789,-3396)(-7.20000,3.60000){5}{\makebox(6.3500,9.5250){\SetFigFont{7}{8.4}{\rmdefault}{\mddefault}{\updefault}.}}
\multiput(6760,-3382)(-7.90000,3.95000){5}{\makebox(6.3500,9.5250){\SetFigFont{7}{8.4}{\rmdefault}{\mddefault}{\updefault}.}}
\multiput(6728,-3367)(-8.62070,3.44828){5}{\makebox(6.3500,9.5250){\SetFigFont{7}{8.4}{\rmdefault}{\mddefault}{\updefault}.}}
\multiput(6694,-3352)(-7.41380,2.96552){6}{\makebox(6.3500,9.5250){\SetFigFont{7}{8.4}{\rmdefault}{\mddefault}{\updefault}.}}
\multiput(6657,-3337)(-7.75862,3.10345){6}{\makebox(6.3500,9.5250){\SetFigFont{7}{8.4}{\rmdefault}{\mddefault}{\updefault}.}}
\multiput(6618,-3322)(-7.93104,3.17242){6}{\makebox(6.3500,9.5250){\SetFigFont{7}{8.4}{\rmdefault}{\mddefault}{\updefault}.}}
\multiput(6578,-3307)(-8.17242,3.26897){6}{\makebox(6.3500,9.5250){\SetFigFont{7}{8.4}{\rmdefault}{\mddefault}{\updefault}.}}
\multiput(6537,-3291)(-8.22000,2.74000){6}{\makebox(6.3500,9.5250){\SetFigFont{7}{8.4}{\rmdefault}{\mddefault}{\updefault}.}}
\multiput(6496,-3277)(-8.04000,2.68000){6}{\makebox(6.3500,9.5250){\SetFigFont{7}{8.4}{\rmdefault}{\mddefault}{\updefault}.}}
\multiput(6456,-3263)(-9.30000,3.10000){5}{\makebox(6.3500,9.5250){\SetFigFont{7}{8.4}{\rmdefault}{\mddefault}{\updefault}.}}
\multiput(6419,-3250)(-8.02500,2.67500){5}{\makebox(6.3500,9.5250){\SetFigFont{7}{8.4}{\rmdefault}{\mddefault}{\updefault}.}}
\multiput(6387,-3239)(-9.30000,3.10000){4}{\makebox(6.3500,9.5250){\SetFigFont{7}{8.4}{\rmdefault}{\mddefault}{\updefault}.}}
\multiput(6359,-3230)(-10.95000,3.65000){3}{\makebox(6.3500,9.5250){\SetFigFont{7}{8.4}{\rmdefault}{\mddefault}{\updefault}.}}
\multiput(6337,-3223)(-8.55000,2.85000){3}{\makebox(6.3500,9.5250){\SetFigFont{7}{8.4}{\rmdefault}{\mddefault}{\updefault}.}}
\multiput(6320,-3217)(-9.90000,3.30000){2}{\makebox(6.3500,9.5250){\SetFigFont{7}{8.4}{\rmdefault}{\mddefault}{\updefault}.}}
{\color[rgb]{0,0,0}\put(7426,-1861){\line( 0,-1){4275}}
}%
\thinlines
{\color[rgb]{0,0,0}\multiput(5626,-4336)(276.25931,331.51117){7}{\line( 5, 6){131.624}}
}%
{\color[rgb]{0,0,0}\multiput(5626,-4336)(328.70979,197.22587){6}{\line( 5, 3){176.304}}
}%
{\color[rgb]{0,0,0}\multiput(5626,-4336)(319.62321,-319.62321){6}{\line( 1,-1){145.384}}
}%
{\color[rgb]{0,0,0}\multiput(5626,-4336)(398.38673,-132.79558){5}{\line( 3,-1){195.053}}
}%
\put(7651,-2198){\makebox(0,0)[lb]{\smash{{\SetFigFont{11}{16.8}{\rmdefault}{\mddefault}{\updefault}{\color[rgb]{0,0,0}$\RP^1$}%
}}}}
{\color[rgb]{0,0,0}\put(5626,-5461){\vector( 0, 1){3600}}
}%
\end{picture}%
\end{center}

Part (ii). The proof will be given in Section~\ref{HiFr}, see Corollary~\ref{equi}.
\proofend

Clearly, $SL(2,\R)$-equivalent curves $\Gamma$ give rise to the same function $F(x,y)$.

\begin{example} \label{sinus}
{\rm If $\Gamma$ is an arc length parameterized unit circle, we obtain the continuous frieze pattern $F(x,y)=\sin(y-x)$ with $T=\pi$.
}
\end{example}

\begin{remark} \label{gen}
{\rm If one does not care about boundary conditions, then solutions to (\ref{Liou}) can be obtained from two curves, $\Gamma$ and $\widetilde\Gamma$, both satisfying (\ref{deter}):
$$
F(x,y)=[\Gamma(x),\widetilde\Gamma(y)].
$$
For example, if 
$
\Gamma(x)=(x,-1),\, \widetilde\Gamma(y)=(1,y),
$
we obtain $F(x,y)=1+xy$. 
A more general solution of this kind is $F(x,y)=(xy)^t + (xy)^{1-t}$ for any real $t$.
}
\end{remark}

\subsection{Hill's equations and projective curves from continuous friezes}\label{HiFr}

Let us show that every closed continuous frieze pattern can be obtained from a projective curve.
In other words, we will show that the construction of Section~\ref{ConSSec} is universal.

\begin{lemma} 
\label{eqfromf}
Let $F(x,y)$ be a closed continuous frieze pattern. 
Then $F$, as a function of $x$ only, or as a function of $y$ only, satisfies the same Hill equation
$$
F_{xx}(x,y)=k(x) F(x,y),\qquad 
F_{yy}(x,y)=k(y) F(x,y)
$$
with $T$-periodic potential $k$ and monodromy $-\Id$.
\end{lemma}

\proof Differentiate (\ref{Liou}) to obtain
$
F F_{xxy}=F_{y} F_{xx}.
$
Thus 
$$
F_{xx}(x,y)=k(x,y) F(x,y),\qquad F_{xxy}(x,y)=k(x,y) F_{y}(x,y)
$$
for some function $k(x,y)$. 
Differentiate the first of these equations with respect to $y$ to obtain \allowbreak $k_{y}(x,y)=0$. 
Hence $k$ depends on $x$ only:
$$
F_{xx}(x,y)=k(x) F(x,y),\qquad
F_{yy}(x,y)=m(y) F(x,y),
$$
where the second equation is obtained similarly to the first one.

To prove that $k=m$, differentiate the second equality in (\ref{bdry}) to obtain
$$
F_{xy}(x,x)+F_{yy}(x,x)=0,\qquad F_{xx}(x,x)+F_{xy}(x,x)=0,
$$
and hence $F_{xx}(x,x)=F_{yy}(x,x)$. This implies that $k(x)=m(x)$.

The third equality in (\ref{bdry}) implies that the monodromy of the Hill equation is $-\Id$, and it follows that $k$ is $T$-periodic.
\proofend

\begin{corollary} 
\label{equi}
The above constructions provide a bijection between projective equivalence classes of 
$T$-periodic parameterizations of $\RP^1$ and closed continuous frieze patterns.
\end{corollary}

\proof
First of all, a projective equivalence class of a curve in $\RP^1$ is the same as Hill's equation; 
see Proposition~\ref{KirProp}. 
Start with a curve $\Gamma(x)$ satisfying Hill's equation. 
Then $F(x,y)=[\Gamma(x),\Gamma(y)]$, and one has:
$$
F_{xx}(x,y)=[\Gamma''(x),\Gamma(y)]=k(x) F(x,y),
$$ 
that is, one recovers the Hill equation from the respective continuous frieze, and thus the 
$\SL_2(\R)$-equivalent class of the curve $\Gamma$.

Conversely, let us show that the function $k$ uniquely determines a continuous frieze. 
Indeed, $F(x,y)$, as a function of~$x$, is the solution of Hill's equation $f''(x)=k(x) f(x)$ 
with the initial conditions $f(y)=0, f'(y)=-1$.
\proofend

This completes the proof of Theorem~\ref{curve}.

\section{The symplectic structure} \label{SSect}

In the previous section, we showed that the space $\Diff_+(S^1)/\PSL_2(\R)$ is
a continuous limit of the space of Coxeter's frieze patterns.
In this section, we will compare the symplectic structures on both spaces
and complete the proof of Theorem~\ref{MainOne}.

\subsection{Kirillov's symplectic structure} \label{Kiss}
The space $\Diff_+(S^1)/\PSL_2(\R)$ is a coadjoint orbit of the Virasoro algebra
and therefore it has the canonical Kirillov symplectic $2$-form.
Let us give here several equivalent expressions of this $2$-form;
see~\cite{K1,K2,OT} for more details.

The first expression is just the definition of the Kirillov symplectic form.
Given a Hill equation (\ref{Hill}) with potential $k(x)$,
it is identified with an element of the dual space to the Virasoro algebra
as follows.
Let $(X(x)\frac{d}{dx},\alpha)$ be any element of the Virasoro algebra,
then 
$$
\left\langle
\left(k(x),c\right),\,\Big(X(x)\frac{d}{dx},\alpha\Big)\right\rangle :=
\int_0^Tk(x)X(x)\,dx+c\alpha.
$$
The coadjoint action of the Virasoro algebra is given by:
\begin{equation}
\label{smalladStar}
ad^*_{X(x)\frac{d}{dx}}\left(k(x),c\right)=
\left(X(x)k'(x)+2X'(x)k(x)+c\,X'''(x),0\right),
\end{equation}
which is nothing else but the infinitesimal version of~(\ref{AdStar}).
Every tangent vector to the coadjoint orbit of the Virasoro algebra
through the point $\left(k(x),c\right)$ is obtained by the coadjoint action
of some vector field.
Let
$$
\xi=ad^*_{X(x)\frac{d}{dx}}\left(k(x),c\right)
\qquad\hbox{and}\qquad 
\eta=ad^*_{Y(x)\frac{d}{dx}}\left(k(x),c\right),
$$ 
then by definition,
\begin{equation}
\label{OmeK}
\begin{array}{rcl}
\om_K(\xi,\eta)&=&
\displaystyle
\int_0^T k(x)\left(X(x)Y'(x)-X'(x)Y(x)\right)dx-
c\,\int_0^TX'(x)Y''(x)\,dx\\[10pt]
&=&
\displaystyle
-\int_0^T\left(X(x)k'(x)+2X'(x)k(x)+c\,X'''(x)\right)Y(x)\,dx.
\end{array}
\end{equation}
Note that the last term on the first row of~(\ref{OmeK}) is the famous Gelfand-Fuchs cocycle.
This formula is independent of the choice of the coordinate $x$.


Assume that the Hill equation is non-oscillating, 
let us give an equivalent formula for Kirillov's symplectic structure,
in terms of projective curves.
Let $\gamma(x)$ be a projective curve and $\Gamma(x)$ its lift to $\R^2$ satisfying~(\ref{deter}).
For an arbitrary choice of linear coordinates 
$\Gamma(x)=(\Gamma_1(x),\Gamma_2(x))$, consider the function 
$f(x)=\frac{\Gamma_1(x)}{\Gamma_2(x)}$ which determines the curve~$\gamma$.
Let, as above, $\xi$ and $\eta$ be two tangent vectors.
Viewed as variations of~$f(x)$, these tangent vectors are expressed as $T$-periodic functions,
$\xi(x)$ and $\eta(x)$.

\begin{lemma}
\label{CurKir}
One has
\begin{equation}
\label{OmeCur}
\om_K(\xi,\eta)=-c\,\int_0^T\frac{\xi'(x)\eta''(x)-\xi''(x)\eta'(x)}{\left(f'(x)\right)^2}dx.
\end{equation}
This formula does not depend on the choice of the parameter $x$.
\end{lemma}

\proof
The coadjoint action~(\ref{smalladStar}) of the Virasoro algebra,
written in terms of projective curves, reads simply as
$ad^*_{X(x)\frac{d}{dx}}\left(f(x)\right)=X(x)f'(x).$
Therefore,
$$
\xi(x)=\frac{X(x)}{f'(x)},
$$
and similarly for $\eta$. 
Substutute these expressions to~(\ref{OmeK}), and integrate by parts
taking into account $k=cS(f)$, to obtain the result.
It is then easy to check that changing the parameter leaves the formula intact.
\proofend

\subsection{The cluster symplectic structure} \label{Cluss}
The space of closed Coxeter's frieze patterns of width $w=n-3$ is an
algebraic variety of dimension $w$.
It has a structure of cluster manifold and therefore has a canonical 
closed $2$-form, i.e., a (pre)symplectic form; see~\cite{GSV} for a general theory.
Let us give the explicit expression of this $2$-form.

\begin{definition}
{\rm
Given a coordinate system $(a_1,a_2,\ldots,a_w)$
associated to an South-East diagonal~(\ref{CanDiag}),
the canonical cluster symplectic form on the space of friezes
is defined by the formula
\begin{equation}
\label{CanSym}
\om=\sum_{1\leq{}i\leq{}w-1}
\frac{da_i\wedge{}da_{i+1}}{a_ia_{i+1}}.
\end{equation}
}
\end{definition}

Consider now an arbitrary zigzag~(\ref{ZZCo}) coordinates
$(a'_1,a'_2,\ldots,a'_w)$.
Define the following {\it a-priori} different $2$-form
$$
\om'=\sum_{1\leq{}i\leq{}w-1}
(-1)^{\varepsilon_i}\,\frac{da'_i\wedge{}da'_{i+1}}{a'_ia'_{i+1}},
$$
where $\varepsilon_i=0$ if $(a'_i,a'_{i+1})$ belongs to a South-East diagonal,
and $\varepsilon_i=1$ if $(a'_i,a'_{i+1})$ belongs to a South-West diagonal.

\begin{proposition}
\label{Canon}
One has $\om=\om'$, for any zigzag coordinate system.
\end{proposition}

\proof It suffices to examine how the 2-form changes under an elementary transformation of a zigzag $\dots bac \dots \mapsto \dots bcd \dots$ in (\ref{diamond1}). In this case, the difference $\om-\om'$ belongs to the ideal generated by 
the differential of the defining identity of frieze pattern (\ref{RulEq}).
\proofend

The form (\ref{CanSym}) is symplectic, i.e., it is non-degenerate,
if and only if $w$ is even (that is, $n$ is odd).
Otherwise, the form $\om$ has a kernel of dimension $1$.

\subsection{Continuous limit of the cluster symplectic form} \label{ClimOm}
We will now apply the procedure of taking the continuous limit
from Section~\ref{ConLimS} to obtain the continuous limit
of the symplectic form~(\ref{CanSym}).

The symplectic form $\om$ written in geometric terms (see Remark~\ref{ImpRem}) is as follows:
$$
\om=\sum_{1\leq{}i\leq{}w-1}
\frac{d\left[V_{n-1},V_i\right]\wedge{}d\left[V_{n-1},V_{i+1}\right]}
{\left[V_{n-1},V_i\right]\left[V_{n-1},V_{i+1}\right]}.
$$
Let $\xi,\eta$ be two tangent vectors to the space of friezes.
Each of them can be represented by $n$ vectors in $\R^2$, that is,
$
\xi=\left(
\xi_0,\ldots,\xi_{n-1}
\right),
$
such that
$$
\left[V_{i},\xi_{i+1}\right]+\left[\xi_i,V_{i+1}\right]=0,
$$
since $\left[V_{i},V_{i+1}\right]\equiv1$.
We obtain
$$
\om(\xi,\eta)=\sum_{1\leq{}i\leq{}w-1}
\frac{\left[V_{n-1},\xi_i\right]\left[V_{n-1},\eta_{i+1}\right]-
\left[V_{n-1},\xi_{i+1}\right]\left[V_{n-1},\eta_{i}\right]}
{\left[V_{n-1},V_i\right]\left[V_{n-1},V_{i+1}\right]}.
$$

The continuous limit of the $n$-gon $(V_0,\ldots,V_{n-1})$ is 
a curve $\Gamma(x)=\left(\Gamma_1(x),\Gamma_2(x)\right)$.
A tangent vector is represented by a curve
$\xi(x)=\left(\xi_1(x),\xi_2(x)\right)$ such that 
$$
\left[\Gamma,\xi'\right]+\left[\xi,\Gamma'\right]=0.
$$
The continuous limit of the above sum is then the following integral:
$$
\om(\xi,\eta)=\int_0^\pi\frac{\xi_2(x)\eta'_2(x)-\xi'_2(x)\eta_2(x)}{\Gamma_2(x)^2}dx.
$$
Let us show that this expression coincides with~(\ref{OmeCur})
up to the multiple $-\frac{1}{4c}$.

A projective curve can be thought of as a function $f(x)=\frac{\Gamma_1(x)}{\Gamma_2(x)}$.
The affine lift $(f(x),1)$ does not satisfy the equality (\ref{deter}) but the rescaling
\begin{equation}
\label{GamLift}
\Gamma(x)=\left(f(x)f'(x)^{-1/2}, f'(x)^{-1/2}\right)
\end{equation}
does.
Let $\xi(x)$ be a tangent vector, i.e., a variation of the function $f(x)$.
Lifted to a tangent vector on curves $\xi(x)=(\xi_1(x),\xi_2(x))$ it then reads  as
$$
\left(\xi_1(x),\xi_2(x)\right)=
\left(\xi(x)f'(x)^{-1/2}-\frac{1}{2}\xi'(x)f'(x)^{-3/2},\,-\frac{1}{2}\xi'(x)f'(x)^{-3/2}\right).
$$
One readily obtains $\om(\xi,\eta)=-\frac{1}{4c}\om_K(\xi,\eta)$.

We have proved that the continuous limit of the cluster (pre)symplectic form
on the space of Coxeter's friezes is (up to a multiple) the Kirillov
symplectic form on $\Diff_+(S^1)/\PSL_2(\R)$.
Theorem~\ref{MainOne} is proved.

\section*{Appendix: relation to metrics of constant curvature}

Let us give yet another geometric interpretation of continuous frieze patterns.
Using (\ref{GamLift}) and (\ref{Constr}), we obtain a general form of solution:
\begin{equation} 
\label{genform}
F(x,y)=\frac{f(y)-f(x)}{\sqrt{f'(x)f'(y)}}
\end{equation}
satisfying the first two boundary conditions (\ref{bdry}). 

\begin{example}
{\rm
Chosing $f(x)=\tan x$ yields Example \ref{sinus}. 
If $f(x)=x$, we obtain a linear solution $F(x,y)=y-x$. }
\end{example}

Now we describe the relation of continuous frieze patterns with conformal metrics 
of constant curvature in dimension $2$.

\begin{lemma}
\label{curv}
The conformal metric $-4F^{-2}(z,\bar z) dz d\bar z$ has constant curvature $-1$ if and only if the function $F$ satisfies equation
(\ref{Liou}).
\end{lemma}

\proof The curvature of the  metric $g=g(z,\bar z) dz d\bar z$ equals 
$$
-\frac{2}{g} \frac{\partial^2\ln g}{\partial z \partial \bar z},
$$
see, e.g., \cite{DFN}. Substituting $g=-4F^{-2}$ yields the result.
\proofend

\begin{example}
{\rm
The Poincar\'e  half plane  metric
\begin{equation}
\label{uphalf}
g=-4\frac{dz d\bar z}{(z-\bar z)^2}
\end{equation}
gives  $F(x,y)=y-x$. More generally, start with (\ref{uphalf}) and change the variable $z=f(w)$. This yields the formula
$$
g=-4\frac{f'(w)f'(\bar w)dw d\bar w}{(f(w)-f(\bar w))^2},
$$
and by Lemma \ref{curv},  the solution (\ref{genform}).}
\end{example}

\begin{remark} \label{hyperb}
{\rm One can also consider a Lorentz metric $4F^{-2}(x,y) dx dy$. Then equation
(\ref{Liou}) is equivalent to this metric having curvature 1.  In particular, such conformal Lorentz metrics of constant curvature on the hyperboloid are studied in \cite{KS}; see also \cite{DG,DO}.
}
\end{remark}

\bigskip
{\bf Acknowledgments}.  
We are grateful to Boris Khesin and Sophie Morier-Genoud for enlightening discussions.
V.~O. was partially supported by the PICS05974 ``PENTAFRIZ'' of CNRS. 
S.T. was partially supported by the NSF grant DMS-1105442.


\begin{thebibliography}{99}

\bibitem{ARS} 
I. Assem, C. Reutenauer, D. Smith, 
{\it Friezes.} Adv. Math. {\bf 225} (2010), 3134--3165.

\bibitem{CaCh}
P. Caldero, F. Chapoton, 
{\it Cluster algebras as Hall algebras of quiver representations},
Comment. Math. Helv.  {\bf 81}  (2006),   595--616.

\bibitem{CoCo}
J. H. Conway, H. S. M. Coxeter, 
{\it Triangulated polygons and frieze patterns},
Math. Gaz. {\bf 57} (1973), 87--94 and 175--183.

\bibitem{Cox} 
H. S. M. Coxeter. 
{\it Frieze patterns}.  Acta Arith.  {\bf 18}  (1971), 297--310.

\bibitem{DFN} B. Dubrovin, A. Fomenko, S. Novikov. 
{\it Modern geometry -- methods and applications.} Part I. 
The geometry of surfaces, transformation groups, and fields. 
Second edition. Springer-Verlag, New York, 1992.

\bibitem{DG}
C. Duval, L. Guieu. 
{\it The Virasoro group and Lorentzian surfaces: the hyperboloid of one sheet.} 
J. Geom. Phys. {\bf 33} (2000),  103--127.

\bibitem{DO} 
C. Duval, V. Ovsienko. {\it Lorentz world lines and the Schwarzian derivative.} 
Funct. Anal. Appl. {\bf 34} (2000), 135--137.

\bibitem{FT} 
L.D. Faddeev, L.A. Takhtajan, {\it Liouville model on the lattice.} Lect. Notes in Phys. {\bf 246} (1986),
166--179.

\bibitem{FRS} 
E. Frenkel, N. Reshetikhin, M.A Semenov-Tian-Shansky, 
{\it Drinfeld-Sokolov reduction for difference operators and deformations of W-algebras. I. 
The case of Virasoro algebra.} Commun. Math. Phys. {\bf 192} (1998), 605--629.

\bibitem{GSV}
M. Gekhtman, M. Shapiro, A. Vainshtein, 
{ Cluster algebras and Poisson geometry.} 
Amer. Math. Soc., Providence, RI, 2010.

\bibitem{K1} 
A.A. Kirillov, 
{\it Infinite-dimensional Lie groups: their orbits, invariants and representations. The geometry of moments.}
Lecture Notes in Math., 970, Springer, Berlin, 1982, 101--123.

\bibitem{K2} 
A.A. Kirillov, 
{\it K\"ahler structure on the $K$-orbits of a group of diffeomorphisms of the circle.} 
Funktsional. Anal. Appl. {\bf 21}:2 (1987), 42--45.

\bibitem{KS}
B. Kostant, S. Sternberg.
{\it The Schwarzian derivative and the conformal geometry of the Lorentz hyperboloid},
 Math. Phys. Stud. {\bf 10}, Kluwer Acad. Publ., Dordrecht, 1988, 113--125.
 
 \bibitem{MS} 
I. Marshall, M. Semenov-Tian-Shansky, 
{\it Poisson groups and differential Galois theory of Schroedinger equation on the circle.} 
Comm. Math. Phys. {\bf 284} (2008), 537--552.
 
 \bibitem{MGOT} 
S. Morier-Genoud, V. Ovsienko, S. Tabachnikov. 
{\it 2-Frieze patterns and the cluster structure of the space of polygons.} Ann. Inst. Fourier {\bf 62} (2012), 937--987.

\bibitem{MGOST} 
S. Morier-Genoud, V. Ovsienko, R. Schwartz, S. Tabachnikov. 
{\it Linear difference equations, frieze patterns and combinatorial Gale transform},
 arXiv:1309.3880.
 
 \bibitem{OST} 
V. Ovsienko, R. Schwartz, S. Tabachnikov, 
{\it Liouville-Arnold integrability of the pentagram map on closed polygons},
Duke Math. J., {\bf 162} (2013), 2149--2196.

\bibitem{OT} 
V. Ovsienko, S. Tabachnikov.
{\it Projective differential geometry, old and new:
from the Schwarzian derivative to the cohomology of diffeomorphism groups.}  Cambridge University Press, Cambridge, 2005.

\bibitem{What} 
V. Ovsienko, S. Tabachnikov.
{\it What is $\dots$the Schwarzian derivative?} 
Notices Amer. Math. Soc. {\bf 56}:1 (2009), 34--36.

\end{thebibliography}
\end{document}